\def\F{{\bf F}}
\def\ra{\rightarrow}
\def\CH{{\rm CH}}
\def\Z{{\bf Z}}
\def\P{{\bf P}}
\def\Q{{\bf Q}}
 \def\et{{\acute et}} 
 \def\ind{{\rm ind}}

{\bf Zero-cycles and rational points on some   surfaces over a global function field}

\medskip

J.-L. Colliot-Th\'el\`ene and Sir Peter Swinnerton-Dyer

\medskip

{\it R\'esum\'e}  Soit $\F$ un corps fini de caract\'eristique $p$. Pour une surface lisse sur $\F(t)$ 
d\'efinie par une  \'equation $f+tg=0$,
o\`u $f $ et $g$ sont deux formes de degr\'e $d$ sur $\F$ en 4 variables, avec $d$ premier  \`a $p$, nous montrons que
l'obstruction de Brauer-Manin au principe de Hasse pour les z\'ero-cycles de degr\'e 1
 est la seule obstruction. 
  Pour $d=3$ (surfaces cubiques), 
 on  en d\'eduit le m\^{e}me \'enonc\'e pour les points rationnels.
 
 {\it Summary} Let $\F$ be a finite field of characteristic $p$. We consider smooth surfaces over $\F(t)$ defined by
an equation $f+tg=0$, where $f $ and $g$  are forms of degree $d$ in 4 variables  with coefficients in $\F$, with $d$
prime to $p$.
We prove : For such surfaces over $\F(t)$, the Brauer-Manin obstruction to the existence
of a zero-cycle of degree one  is the only obstruction. For $ d=3 $ (cubic surfaces), this leads to
the same result  for rational points.

\bigskip

{\bf \S 1.  Introduction}

\medskip

Let $k$ be a global field.

\medskip

Study of the case of curves (Cassels, Tate)  and of the case of 
rational surfaces (Colliot-Th\'el\`ene et  Sansuc  [CT/S81], 
where a more precise conjecture is made for
rational surfaces) has led to the following conjecture
for {\it zero-cycles} on arbitrary   varieties over global fields
(Kato/Saito [K/S86], Saito [S89], Colliot-Th\'el\`ene [CT93], [CT99]).

 \medskip

{\bf Conjecture}  {\sl Let   $X$ be a smooth, projective, geometrically integral
variety over a global field
 $k$. If there exists a family
$\{z_v\}_{v \in \Omega}$ of local zero-cycles of degree 1
(here $v$ runs through the set $\Omega$ of places of $k$)
such that for all
$A \in {\rm Br}(X)$,  
$$\sum_{v \in \Omega} {\rm inv}_v(A(z_v))=0 \in \Q/\Z$$  holds, then there exists
a zero-cycle of degree 1 on $X$.
 In other words, the Brauer--Manin obstruction to the existence of a
 zero-cycle of degree 1 on $X$ is the only obstruction.}

 \medskip
 
 Over number fields, this conjecture has been established in special cases
in work of (alphabetical order, and various combinations) 
Colliot-Th\'el\`ene, Frossard, Salberger,  Sansuc, Skorobogatov,
Swinnerton-Dyer, Wittenberg (see the introduction of [W10]).
None of these results applies to smooth surfaces of degree at least 3
in 3-dimensional projective space -- for $d\geq 5$ these surfaces are of general type.
In section 2,   we establish the conjecture in the    special case of a global field $k=\F(t)$ purely transcendental over a finite field $\F$ and of smooth surfaces $X \subset \P^3_{k}$ defined by an equation $f+tg=0$, where $f$ and $g$ are two   forms of arbitrary  degree $d$ over the field $\F$.

\bigskip

 According to a conjecture of Colliot-Th\'el\`ene
and Sansuc ([CT/S80]), the Brauer--Manin obstruction to the existence of a {\it rational
point} on a smooth, geometrically rational surface defined over a global field should be the only obstruction. Such should in particular be the case for smooth cubic surfaces in 3-dimensional
projective space $\P^3_{k}$.    
In section 3, we establish  the conjecture in the   special case of a global field $k=\F(t)$ purely transcendental over a finite field $\F$ and of smooth cubic surfaces $X \subset \P^3_{k}$ defined by an equation $f+tg=0$, where $f$ and $g$ are two cubic forms over the field $\F$. Simple though they be, such surfaces may fail to obey the Hasse principle.

\bigskip

{\bf  \S 2.  Zero-cycles of degree 1 on surfaces of arbitrary degree}

\medskip

The following theorem is due to S. Saito [S89]. It says that if a strong integral form
of the Tate conjecture on 1-dimensional cycles is true, then the above conjecture holds,
at least if we stay away from the characteristic of the field.
For an alternate proof of Theorem 1, see    [CT99], Prop. 3.2.

\medskip
\vfill\eject

{\bf Theorem 1} (Saito)  {\sl Let  $\F$ be a finite field and 
$C/\F$  a smooth, projective, geometrically integral curve over $\F$.
Let $k=\F(C)$ be its function field.
Let 
${\cal X}$ be a smooth, projective, geometrically integral
$\F$-variety of dimension
$n$ and $f : {\cal X} \ra C$  a faithfully flat map whose generic fibre
$X/k$ is smooth and geometrically integral.

Assume :

(1)  For each prime $l \neq {\rm char} (\F)$,
the cycle map
$$T_{X} : \CH^{n-1}({\cal X}) {\otimes}\Z_l \ra H^{2n-2}_{\et}({\cal X},\Z_l(n-1))$$
from the Chow group of dimension 1 cycles on $\cal X$ to \'etale cohomology
 is onto.
 
 (2) There exists a family
$\{z_v\}_{v \in \Omega}$ of local zero-cycles of degree 1
(here $v$ runs through the set $\Omega$ of places of $k$)
such that for all
$A \in {\rm Br}(X)$,  
$$\sum_{v \in \Omega} {\rm inv}_v(A(z_v))=0 \in \Q/\Z.$$
 Then there exists a zero-cycle on $X$ of degree 
  a power of ${\rm char} (\F)$.}
  
  \medskip

In this statement, $A(z_{v})$ is the element of the Brauer group of the local field $k_{v}$
obtained by evaluation of $A$ on the zero-cycle $z_{v}$. The map ${\rm inv}_{v} : {\rm Br}(k_{v}) \to \Q/\Z$
is the local invariant of class field theory.
 
\medskip

Here is one case where assumption (1) in the previous theorem is fulfilled.

\medskip

{\bf Theorem 2} {\sl Let $\F$ be a finite field and $l$ a prime, $l \neq {\rm char}(\F)$.
For a smooth, projective, geometrically integral threefold ${\cal X}$
over  $\F$ which is birational to $\P^{3}_{F}$, the cycle map
$T_{\cal X} :  CH^2({\cal X}) {\otimes}\Z_l \ra H^{4}_{\et}({\cal X},\Z_l(2))$
is onto.}

\medskip

{\it Proof }
If ${\cal X} = \P^{3}_{\F}$, then $CH^2({\cal X})=\Z$
and one easily checks that the  cycle map
$$ T_{\cal X} : CH^2({\cal X}) {\otimes}\Z_l \ra H^{4}_{\et}({\cal X},\Z_l(2))$$
is simply the identity map $\Z_{l}=\Z_{l}$.
 Using the standard formulas for the computation of Chow groups
and of cohomology for a blow-up along a smooth projective 
subvariety, as well as the vanishing
 of  Brauer groups of smooth projective 
curves over a finite field, one shows :
For ${\cal X}$ a smooth projective threefold, the cokernel of the above cycle map $T_{\cal X}$  is 
invariant under blow-up of  smooth projective subvarieties on ${\cal X}$.

By a result of Abhyankar ([Abh66], Thm. 9.1.6), 
there exists a smooth projective variety
${\cal X}'$ which is obtained from $\P^3_{\F}$ 
by a sequence of blow-ups along smooth projective
$\F$-subvarieties, and which is equipped with  a birational $\F$-morphism $p : {\cal X}' \to {\cal X}$.

There are push-forward maps $\pi_{*}$ and pull-back maps $\pi^*$
both for Chow groups and for \'etale cohomology, and for the birational
map $\pi$ we have $\pi_{*} \circ \pi^*= {\rm id}$. Moreover these maps
are compatible with the cycle class map. 
Thus the cokernel of  $T_{\cal X}$ is a subgroup of the cokernel of
$T_{{\cal X}'}$, hence is zero. QED

\bigskip

 Combining Theorems 1 and 2, we get :
 
 \medskip

{\bf Theorem 3}  {\sl Let  $\F$ be a finite field and 
$C/\F$  a smooth, projective, geometrically integral curve over $\F$.
Let $k=\F(C)$ be its function field.
Let 
${\cal X}$ be a smooth, projective, geometrically integral
$\F$-variety of dimension
$n$ and $f : {\cal X} \ra C$  a faithfully flat map whose generic fibre
$X/k$ is smooth and geometrically integral.

Assume :

(1)   ${\rm dim} \ {\cal X}= 3$  
and ${\cal X}$ is $\F$-rational;

 (2) there exists a family
$\{z_v\}_{v \in \Omega}$ of local zero-cycles of degree 1
(here $v$ runs through the set $\Omega$ of places of $k$)
such that for all
$A \in {\rm Br}(X)$,  
$$\sum_{v \in \Omega} {\rm inv}_v(A(z_v))=0 \in \Q/\Z.$$
 
Then there exists a zero-cycle on $X$ of degree 
a power of  ${\rm char} (\F)$.}

\medskip
\vfill\eject

We may now prove the main result of this section.

\medskip

{\bf Theorem  4} {\sl Let  $\F$ be a finite field, let $f,g$ be two nonproportional homogeneous forms in 4 variables,
of degree $d$ prime to the characteristic of $\F$. 
 Let $k=\F(t)$.
Suppose the $k$-surface $X \subset \P^3_{k}$ defined by $f+tg=0$ is smooth.
 If there is no Brauer--Manin obstruction to the Hasse principle for zero-cycles
of degree 1 on $X$, then

(i) there exists a zero-cycle of degree 1 on the $k$-surface $X$;

(ii) there exists a zero-cycle of degree 1 
on the $\F$-curve $\Gamma$ defined by $f=g=0$ in $\P^3_{\F} $.}

\medskip

{\it Proof}  Let ${\cal X}_{1} \subset \P^3_{\F}\times_{F}\P^1_{\F}$ be the schematic closure
of $X \subset  \P^3_{\F(t)}$.  The $\F$-variety  ${\cal X}_{1}$  has an affine birational model with equation
 $\phi(x,y,z) +t \psi(x,y,z)=0$, hence $t$ is determined
by $x,y,z$, thus ${\cal X}$ is $\F$-birational to $\P^3_{\F}$.
Since ${\cal X}_{1}$ admits a smooth projective model over $\F$,
 a result  of Cossart  ([Co92], Th\'eor\`eme, p.~115)  shows that there exists  a smooth projective threefold ${\cal X}/\F$  and an $\F$-birational morphism  ${\cal X} \to {\cal X}_{1}$ which is an isomorphism over the smooth locus
of ${\cal X}_{1}$, hence in particular
which induces an isomorphism  over $Spec \  \F(t) \subset \P^1_{\F}$.
That is, the generic fibre of ${\cal X} \to \P^1_{\F}$ is $k$-isomorphic to $X/k$.

It remains to combine Theorem 1 and Theorem 2 to prove (i).
Statement (ii) follows from (i) as a special application
of a result of Colliot-Th\'el\`ene and Levine ([CT/L09], Th\'eor\`eme 1,  p.~217). 
 QED

\medskip

 {\bf Remark}
Theorem 4  is of interest only in the case where
the $\F$-curve $\Gamma$ does not contain a geometrically
integral component. Otherwise the two statements immediately follow
from  the Weil estimates for  the number of points on geometrically integral curves.
These estimates actually provide more : they show that
if there exists such a component, then
on any field extension  $\F'$ of $\F$ of high enough degree,
there exists an $\F'$-point on $\Gamma$, hence for any such field there exists
an $\F'(t)$-point on the $\F(t)$-surface $X$.

\medskip

{\bf Remark}
One could try to circumvent the cohomological machinery, i.e. Theorems 1
and  2.
For this, in each of the special cases where there are zero-cycles of degree 1
everywhere locally on $X$  but there is no zero-cycle of degree one on
the curve $\Gamma$, one should:

(i) Check that the Brauer group is not trivial, find generators.

(ii) Check that there is a Brauer--Manin obstruction.

Already when the common degree of $f$ and $g$ is 3, 
which we shall now more particularly examine,
this seems no easy enterprise.

\bigskip

{\bf \S 3. Rational points on cubic surfaces}

\medskip

The proof of the  following result is independent of the previous results.

\medskip

{\bf Theorem 5} 
{\sl Let  $\F$ be a finite field, let $f,g$ be two nonproportional  cubic forms over $\F$
 in 4 variables.
Assume the characteristic of $\F$ is not 3. 
 Let $k=\F(t)$.
Suppose the $k$-surface $X \subset \P^3_{k}$ defined by $f+tg=0$ is smooth.
Let  $\Gamma \subset \P^3_{\F}$ be the complete intersection curve defined
by $f=g=0$.
The following conditions are equivalent :

(i) There exists a $k$-rational point on the $k$-variety $X$.

(ii) There exists a zero-cycle of degree 1 on the $k$-variety $X$.

(iii)  
There exists a zero-cycle of degree 1 on the $\F$-curve $\Gamma$.

(iv)  There exists
a closed point of degree prime to 3 on the $\F$-curve $\Gamma$.

(v)  There exists  a closed point of degree a power of 2 on the $\F$-curve $\Gamma$.}
 
 \medskip
 
 {\it Proof}
  That (i) implies (ii)  is trivial.
  That (ii) implies (iii) is a special case of [CT/L09].
  Statements (iii) and (iv) are equivalent, since $\Gamma$ is a curve
  of degree 9.
  If (v) holds, then $\Gamma$ has a point in a tower of quadratic extensions of $\F$,
  hence the cubic surface $X$ has a point in a tower of quadratic extensions of $k$.
  An extremely well known argument shows that if a cubic surface over a field
  has a point in a separable quadratic extension of that field, then it has a rational point :
 the line joining two conjugate points is defined over the ground field, either it is entirely
 contained in the cubic surface or it meets it in a third, rational point. Iterating
 this remark, we see that $X$ has a rational point, i.e. (i) holds.
  
 \medskip
  
 Let us prove that (iii) implies (v). 
 To prove this, one may replace $\F$
 by its maximal multiquadratic extension $F$, which we now do.
 For an odd integer $n$, we let $F_{n}/F$ be the unique  field extension of $F$
 of degree $n$.
 
 For $Z/L$ a variety over a field $L$, the index $\ind(Z)=\ind(Z/L)$ is
 the gcd of the $L$-degrees of closed points on $Z$.
 The index of an $L$-variety is equal to the index of its reduced $L$-subvariety.
The index of an $L$-variety which is a finite union of   $L$-varieties is the gcd of the indices of each of them.

Since $F$ has no quadratic or quartic extension, 
an effective zero-cycle of degree $1,2,4$ contains
an $F$-rational point, and an effective zero-cycle of degree 
3, 6 or 9 either contains an $F$-point or has index a multiple
of 3.

If $\Gamma$ contains a geometrically  integral  component,
then $\Gamma(F) \neq \emptyset$ (Weil estimates, see the remark after Theorem 4).
 
 Suppose $\Gamma$ does not contain a geometrically integral   component.
One then easily checks that the degree 9 curve ${\overline{\Gamma}}$ can break up
only  in one of the following ways :

$9= 3(1+1+1)$
 
$9=2(1+1+1) + (1+1+1)$

$9= (2+2+2) + (1+1+1)$

$9= (1+1+1) + (1+1+1) + (1+1+1)$

$9= (1+ \dots + 1)$   (9 times)

$9 = (3 + 3 +3)$

\medskip

Here $ (a+a+a) $ means the sum of 3 conjugate integral curves of degree $a$
over ${\overline F}$.

An integral curve of degree  2 over ${\overline F}$ is a smooth plane conic, contained in a
well-defined plane.

An integral curve of degree 3 over ${\overline F}$ is either a plane cubic or a smooth twisted cubic.

Let the integral  curve $C \subset \P^3_{F}$ break up as (1+1+1). 
The singular set consists of at most 3 points.
Then either $C(F) \neq \emptyset$
or 3 divides $\ind(C)$. 

Let the integral curve  $C \subset \P^3_{F}$ break up as  (2+2+2).
Each conic is defined over $F_{3}$.
Two  distinct smooth conics on $f=0$  define two distinct   planes,
hence they  intersect in at most 2 geometric   points.
Such points must already be in $F_{3}$.
Thus any closed point in the singular locus of $C$
has degree 1 or 3.
One concludes that either $C(F) \neq \emptyset$  
or  3 divides $\ind(C)$.

Let the integral curve $\Gamma  \subset \P^3_{F}$ break up as $(1+Ê\dots +1)$ (9 times).
 The 9 lines are defined over  $F_{9}$, the degree 9
extension of $F$. So are their intersection points.
This implies that any singular closed point on $\Gamma $
has degree a power of 3. Thus $\Gamma (F) \neq \emptyset$  or  3 divides $\ind(\Gamma )$.

Let the integral  curve $\Gamma \subset \P^3_{F}$ break up as  $(3 + 3 +3)$,
and  assume that this corresponds to a decomposition as three
conjugate plane cubics. Each of these is defined over $F_{3}$.
The intersection number  of two of these cubics  is 3.
The points of intersection
of two such curves are thus defined over $F_{9}$. We conclude that
the singular locus of $\Gamma$ splits over $F_{9}$. This implies
that the degree of any closed point in that locus is a power of 3.
Thus either 
$\Gamma(F) \neq \emptyset$
or 3 divides $\ind(\Gamma)$.

Let the curve $\Gamma \subset \P^3_{F}$ break up as  $(3 + 3 +3)$,
and assume that  $\Gamma$ breaks up as
the sum of three conjugate twisted cubics. 
The curve $\Gamma$ lies on the smooth   cubic surface $X$ over $F(t)$  defined by $f+tg=0$.
 Each twisted curve is defined over $F_{3}$.  Let $\sigma$ be a generator
 of ${\rm Gal}(F_{3}(t)/F(t))$. Write $\Gamma=C+\sigma(C)+\sigma^2C$ on $X_{F_{3}(t)}$.
 Using intersection theory on the smooth surface $X_{F_{3}(t)}$, which is
 invariant under the action of ${\rm Gal}(F_{3}(t)/F(t))$,
 and letting $H$ be the class of a plane section,
 we find
$27 =(3H.3H) = (\Gamma.\Gamma) = 3(C.C)+ 6(C.\sigma(C)).$
 The curve $C$ is a twisted
 cubic, hence a smooth curve of genus zero on the smooth cubic surface $X$,
 whose canonical bundle $K$  is given by $-H$. The formula for the arithmetic
 genus of a curve on a surface, namely $2(p_{a}(C)-1)=(C.C)+(C.K)$
 gives $(C.C)=1$. This implies $(C.\sigma(C))=4$, hence $(\sigma(C).\sigma^2(C))=4$
 and $(\sigma^2(C).C)=4$. 
 Since each of these twisted cubic is defined over $F_{3}$ and since
   $F_{3}$ has no field extension of degree 2 or 4,
this implies that the points of intersection of any two of these  twisted cubics  are defined over $F_{9}$. We conclude that
the singular locus of $\Gamma$ splits over $F_{9}$. This implies
that the degree of any closed point in that locus is a power of 3.
Thus either 
$\Gamma(F) \neq \emptyset$
or 3 divides $\ind(\Gamma)$.

 In all cases we have proved :
Either $\Gamma(F) \neq \emptyset$
or 3 divides $\ind(\Gamma)$.

The assumption $\ind(\Gamma)=1$, made in (iii), now implies
$\Gamma(F) \neq \emptyset$. QED
 
 \bigskip
 
 {\bf Remark}
 If  the order $q$ of the finite field $\F$ is large enough and $f+tg=0$ is soluble in $\F(t)$, a variant of
 the proof for  the equivalence of (iv) and (v) shows that $f+tg=0$  has a solution in
    polynomials of degree at most $5$. This raises an interesting general
    question: are there integers $N(d)$ with the following property? Suppose that
    $G(X_0,\ldots, X_4,t)$  is a polynomial defined over $\F$, homogeneous of degree $3$
    in the $X_i$  and of degree $d$ in $t$; if $ G=0$  is soluble in $\F(t)$, then it has
    a solution in polynomials of degree at most $N(d)$.
    
 \bigskip
 
 We may now prove :
 
 \medskip
 
 {\bf Theorem  6}
{\sl
Let  $\F$ be a finite field, let $f,g$ be two nonproportional cubic forms in 4 variables.
  Assume the characteristic of $\F$ is not $3$.  Let $k=\F(t)$.
Suppose the cubic surface $X \subset \P^3_{k}$ over $k$ defined by $f+tg=0$ is smooth.
 If there is no Brauer--Manin obstruction to the Hasse principle for rational points
  on $X$, then there exists a $k$-rational point on  $X$.}

\medskip

 {\it Proof} Combine Theorem 4  and Theorem 5. QED

\bigskip

{\bf Remark}
Again, it  would be nice to avoid the cohomological machinery, 
i.e. Theorems 1 and 2. 
When $X$ has no rational points over $\F(t)$  but points in all the completions 
of $\F(t)$ 
one should  exhibit an explicit Brauer--Manin obstruction for $X$. 
For this purpose, it would probably be helpful to use [SD93].
  Down to earth computations, which we shall not  insert here,
 have led to the following result.
If a smooth cubic surface $X$ given by $f+tg=0$ is a counterexample to the
Hasse principle over $\F(t)$, then, after replacing $\F$ by
its  maximal pro-2-extension $F$,  the following holds :
When going over to the algebraic closure of $F$, 
the curve $\Gamma$ in the proof of Theorem 5 breaks up
as   
 a sum of 9 conjugate lines, or 
  a sum of three twisted cubics, or 
 a sum of three conjugate conics plus  a sum of three  coplanar   conjugate lines;
when using the word ``conjugate'' we mean that the Galois action is   transitive.
Only in these three cases may we   expect a Brauer--Manin obstruction.

\bigskip

{\it References}

\medskip

[Abh66] S. Abhyankar, {\it Resolution of Singularities of Embedded Surfaces}, Academic Press, New York, 1966.

\medskip

[CT93] J.-L. Colliot-Th\'el\`ene, {\it L'arithm\'etique des z\'ero-cycles},  (expos\'e aux Journ\'ees 
arithm\'etiques de Bordeaux, septembre 93),
Journal de th\'eorie des nombres de Bordeaux {\bf 7} (1995) 51--73.

\medskip

[CT99]  J.-L. Colliot-Th\'el\`ene, {\it Conjectures de type local-global sur l'image de l'application cycle en cohomologie \'etale}, 
 in  Algebraic K-Theory (1997), W. Raskind and C. Weibel  ed., Proceedings of Symposia in Pure Mathematics {\bf 67}, Amer. Math. Soc. (1999)  1--12.

\medskip

[CT/L09]  J.-L. Colliot-Th\'el\`ene et Marc Levine,
 {\it Une version du th\'eor\`eme d'Amer et Brumer pour les z\'ero-cycles},
  in {\it Quadratic forms, linear algebraic groups, and cohomology} 
   (ed. J.-L. Colliot-Th\'el\`ene, R. S.  Garibaldi, R. Sujatha, V. Suresh),
   Developments in mathematics,   Springer-Verlag {\bf 18}  (2010),  215--223.

\medskip

[CT/S80] J.-L. Colliot-Th\'el\`ene et J.-J. Sansuc,  {\it La descente sur les vari\'et\'es rationnelles},  
in Journ\'ees de g\'eom\'etrie alg\'ebrique d'Angers (juillet 1979), \'edit\'e par
A.Beauville, Sijthof and Noordhof (1980) 223--237.

\medskip

[CT/S81]  J.-L. Colliot-Th\'el\`ene et J.-J. Sansuc,   {\it On the Chow groups of certain rational surfaces: a sequel
to a paper of S.Bloch}, Duke Math. J. {\bf 48} (1981) 421--447.

\medskip

[Co92] V. Cossart, {\it Mod\`ele projectif r\'egulier et d\'esingularisation}, Math. Annalen {\bf 293} (1992) 115--121.

\medskip

[K/S86] K. Kato et S. Saito, {\it  Global class field theory of arithmetic schemes}, Contemporary math. {\bf 55}, vol. 1 (1986) 255--331.

\medskip

[S89] S. Saito,  {\it Some observations on motivic cohomology of arithmetic schemes}, Invent. math. {\bf 98} (1989) 371-404.

\medskip

[SD93] Sir Peter Swinnerton-Dyer, {\it The Brauer group of cubic surfaces},
Math. Proc. Camb. Phil. Soc. {\bf 113} (1993) 449--460.

\bigskip

[W10] O. Wittenberg,  {\it Z\'ero-cycles sur  les fibrations au-dessus d'une courbe de genre
quelconque}, \break
http://arxiv.org/abs/1010.1883

\vskip1cm
\vfill\eject

J.-L. Colliot-Th\'el\`ene,

C.N.R.S., 

Math\'ematiques, B\^atiment 425,

Universit\'e Paris-Sud

F-91405 Orsay

France

jlct@math.u-psud.fr

\medskip

Sir Peter Swinnerton-Dyer,

Departement of Pure Mathematics and Mathematical Statistics,

Centre for Mathematical Sciences, 

Wilberforce Road,

 Cambridge CB3 0WA

England

H.P.F.Swinnerton-Dyer@dpmms.cam.ac.uk

\bye